\documentclass[12pt]{article}
\usepackage{amsmath,amssymb,amsbsy,amsfonts,amsthm,latexsym,amsopn,amstext,
                                    amsxtra,euscript,amscd}

\begin{document}




\newfont{\teneufm}{eufm10}
\newfont{\seveneufm}{eufm7}
\newfont{\fiveeufm}{eufm5}
%
%
\newfam\eufmfam
                    \textfont\eufmfam=\teneufm \scriptfont\eufmfam=\seveneufm
                    \scriptscriptfont\eufmfam=\fiveeufm

%
%
\def\frak#1{{\fam\eufmfam\relax#1}}
%


\def\bbbr{{\rm I\!R}} 
\def\bbbc{{\rm I\!C}} 
\def\bbbm{{\rm I\!M}}
\def\bbbn{{\rm I\!N}} 
\def\bbbf{{\rm I\!F}}
\def\bbbh{{\rm I\!H}}
\def\bbbk{{\rm I\!K}}
\def\bbbl{{\rm I\!L}}
\def\bbbp{{\rm I\!P}}
\newcommand{\lcm}{{\rm lcm}}
\def\bbbone{{\mathchoice {\rm 1\mskip-4mu l} {\rm 1\mskip-4mu l}
{\rm 1\mskip-4.5mu l} {\rm 1\mskip-5mu l}}}
\def\bbbc{{\mathchoice {\setbox0=\hbox{$\displaystyle\rm C$}\hbox{\hbox
to0pt{\kern0.4\wd0\vrule height0.9\ht0\hss}\box0}}
{\setbox0=\hbox{$\textstyle\rm C$}\hbox{\hbox
to0pt{\kern0.4\wd0\vrule height0.9\ht0\hss}\box0}}
{\setbox0=\hbox{$\scriptstyle\rm C$}\hbox{\hbox
to0pt{\kern0.4\wd0\vrule height0.9\ht0\hss}\box0}}
{\setbox0=\hbox{$\scriptscriptstyle\rm C$}\hbox{\hbox
to0pt{\kern0.4\wd0\vrule height0.9\ht0\hss}\box0}}}}
\def\bbbq{{\mathchoice {\setbox0=\hbox{$\displaystyle\rm
Q$}\hbox{\raise
0.15\ht0\hbox to0pt{\kern0.4\wd0\vrule height0.8\ht0\hss}\box0}}
{\setbox0=\hbox{$\textstyle\rm Q$}\hbox{\raise
0.15\ht0\hbox to0pt{\kern0.4\wd0\vrule height0.8\ht0\hss}\box0}}
{\setbox0=\hbox{$\scriptstyle\rm Q$}\hbox{\raise
0.15\ht0\hbox to0pt{\kern0.4\wd0\vrule height0.7\ht0\hss}\box0}}
{\setbox0=\hbox{$\scriptscriptstyle\rm Q$}\hbox{\raise
0.15\ht0\hbox to0pt{\kern0.4\wd0\vrule height0.7\ht0\hss}\box0}}}}
\def\bbbt{{\mathchoice {\setbox0=\hbox{$\displaystyle\rm
T$}\hbox{\hbox to0pt{\kern0.3\wd0\vrule height0.9\ht0\hss}\box0}}
{\setbox0=\hbox{$\textstyle\rm T$}\hbox{\hbox
to0pt{\kern0.3\wd0\vrule height0.9\ht0\hss}\box0}}
{\setbox0=\hbox{$\scriptstyle\rm T$}\hbox{\hbox
to0pt{\kern0.3\wd0\vrule height0.9\ht0\hss}\box0}}
{\setbox0=\hbox{$\scriptscriptstyle\rm T$}\hbox{\hbox
to0pt{\kern0.3\wd0\vrule height0.9\ht0\hss}\box0}}}}
\def\bbbs{{\mathchoice
{\setbox0=\hbox{$\displaystyle     \rm S$}\hbox{\raise0.5\ht0\hbox
to0pt{\kern0.35\wd0\vrule height0.45\ht0\hss}\hbox
to0pt{\kern0.55\wd0\vrule height0.5\ht0\hss}\box0}}
{\setbox0=\hbox{$\textstyle        \rm S$}\hbox{\raise0.5\ht0\hbox
to0pt{\kern0.35\wd0\vrule height0.45\ht0\hss}\hbox
to0pt{\kern0.55\wd0\vrule height0.5\ht0\hss}\box0}}
{\setbox0=\hbox{$\scriptstyle      \rm S$}\hbox{\raise0.5\ht0\hbox
to0pt{\kern0.35\wd0\vrule height0.45\ht0\hss}\raise0.05\ht0\hbox
to0pt{\kern0.5\wd0\vrule height0.45\ht0\hss}\box0}}
{\setbox0=\hbox{$\scriptscriptstyle\rm S$}\hbox{\raise0.5\ht0\hbox
to0pt{\kern0.4\wd0\vrule height0.45\ht0\hss}\raise0.05\ht0\hbox
to0pt{\kern0.55\wd0\vrule height0.45\ht0\hss}\box0}}}}
\def\bbbz{{\mathchoice {\hbox{$\sf\textstyle Z\kern-0.4em Z$}}
{\hbox{$\sf\textstyle Z\kern-0.4em Z$}}
{\hbox{$\sf\scriptstyle Z\kern-0.3em Z$}}
{\hbox{$\sf\scriptscriptstyle Z\kern-0.2em Z$}}}}
\def\ts{\thinspace}

\newtheorem{theorem}{Theorem}
\newtheorem{lemma}[theorem]{Lemma}
\newtheorem{claim}[theorem]{Claim}
\newtheorem{cor}[theorem]{Corollary}
\newtheorem{prop}[theorem]{Proposition}
\newtheorem{definition}{Definition}
\newtheorem{question}[theorem]{Open Question}
\newtheorem{remark}[theorem]{Remark}

\def\squareforqed{\hbox{\rlap{$\sqcap$}$\sqcup$}}
\def\qed{\ifmmode\squareforqed\else{\unskip\nobreak\hfil
\penalty50\hskip1em\null\nobreak\hfil\squareforqed
\parfillskip=0pt\finalhyphendemerits=0\endgraf}\fi}

\def\cA{{\mathcal A}}
\def\cB{{\mathcal B}}
\def\cC{{\mathcal C}}
\def\cD{{\mathcal D}}
\def\cE{{\mathcal E}}
\def\cF{{\mathcal F}}
\def\cG{{\mathcal G}}
\def\cH{{\mathcal H}}
\def\cI{{\mathcal I}}
\def\cJ{{\mathcal J}}
\def\cK{{\mathcal K}}
\def\cL{{\mathcal L}}
\def\cM{{\mathcal M}}
\def\cN{{\mathcal N}}
\def\cO{{\mathcal O}}
\def\cP{{\mathcal P}}
\def\cQ{{\mathcal Q}}
\def\cR{{\mathcal R}}
\def\cS{{\mathcal S}}
\def\cT{{\mathcal T}}
\def\cU{{\mathcal U}}
\def\cV{{\mathcal V}}
\def\cW{{\mathcal W}}
\def\cX{{\mathcal X}}
\def\cY{{\mathcal Y}}
\def\cZ{{\mathcal Z}}

\newcommand{\comm}[1]{\marginpar{%
\vskip-\baselineskip 
\raggedright\footnotesize
\itshape\hrule\smallskip#1\par\smallskip\hrule}}





\def\MOV{{\bf{MOV}}}

\hyphenation{re-pub-lished}

\def\ord{{\mathrm{ord}}}
\def\Nm{{\mathrm{Nm}}}
\renewcommand{\vec}[1]{\mathbf{#1}}

\def \F{{\bbbf}}
\def \L{{\bbbl}}
\def \K{{\bbbk}}
\def \Z{{\bbbz}}
\def \N{{\bbbn}}
\def \Q{{\bbbq}}
\def\E{{\mathbf E}}
\def\H{{\mathbf H}}
\def\G{{\mathcal G}}
\def\O{{\mathcal O}}
\def\cS{{\mathcal S}}
\def \R{{\bbbr}}
\def\Fp{\F_p}
\def\Fq{\F_q}
\def \fp{\Fp^*}
\def \fq{\Fq^*}
\def\\{\cr}
\def\({\left(}
\def\){\right)}
\def\fl#1{\left\lfloor#1\right\rfloor}
\def\rf#1{\left\lceil#1\right\rceil}

\def\Zm{\Z_m}
\def\Zt{\Z_t}
\def\Zp{\Z_p}
\def\Um{\cU_m}
\def\Ut{\cU_t}
\def\Up{\cU_p}

\def\ep{{\mathbf{e}}_p}
\def\HH{\cH}

\def\Tr{\mbox{Tr}}

\def \Prob{{\mathrm {}}}

\def\LC{{\cL}_{C,\cF}(Q)}
\def\LCn{{\cL}_{C,\cF}(nG)}
\def\Mrs{\cM_{r,s}\(\F_p\)}

\def\taubar{\overline{\tau}}
\def\Fn{\F_{q^n}}
\def\En{\E(\Fn)}

\def\mand{\qquad \mbox{and} \qquad}


\title{A Note on Stable Quadratic Polynomials over Fields of Characteristic Two}

\author{
{\sc Omran Ahmadi}  \\
{Claude Shannon Institute}\\
{University College Dublin} \\
{Dublin 4, Ireland} \\
{\tt omran.ahmadi@ucd.ie}}

\date{\today}

\maketitle

\begin{abstract}
In this note, first we show that there is no stable quadratic polynomial over finite fields of characteristic two
and then show that there exist stable quadratic polynomials over function fields of characteristic
two.
\end{abstract}

\section{Introduction}
Let $K$ be a field and $K[x]$ be the polynomial ring over $K$. A polynomial $f(x)\in K[x]$ is 
called stable if $f(x),ff(x),fff(x),\ldots,f^{(n)}(x),\ldots$ is a sequence of irreducible 
polynomials in $K[x]$. Recently there has been an interest in the study of stable polynomials
.(see~\cite{Nidal,Ayad-McQ,Ayad-McQ-cor,Boston-Jones-1,Jones-1,Jones-2,Alina-Igor})

In~\cite{Alina-Igor}, Ostafe and Shparlinski studied stable quadratic polynomials over 
finite fields and obtained an upper bound on the length of the critical orbits of stable 
quadratic polynomials over
finite fields with odd characteristic. They also posed the question of estimating the number
of stable quadratic polynomials over finite fields. Here we show that there is no stable quadratic
polynomial over finite fields of characteristic two by proving the following theorem.
\begin{theorem}\label{main}
Let $q=2^m$, and let $f(x)=cx^2+ax+b\in\Fq[x]$ where $\Fq$ denotes the field with 
$q$ elements and $a,b\in\fq$. Then $fff(x)$ cannot be an irreducible polynomial
over $\Fq$.
\end{theorem}

\section{Proof of Theorem 1}

For simplicity we assume that $c=1$ and hence $f(x)$ is a monic polynomial. The same proof
can be used for the non-monic polynomials. We also note that in what follows 
${\mbox{\Tr}}_{|K|||L|}(\cdot)$ denotes the trace map from field $K$ with $|K|$ elements to its
subfield $L$ with $|L|$ elements. Our proof is based on the following two well-known 
lemmas.

\begin{lemma}\label{irre-trace}~\cite[Corollary~3.6]{Menezes}
Let $q=2^m$, and let $f(x)=x^2+ax+b\in\Fq[x]$ where $\Fq$ denotes the field with 
$q$ elements and $a,b\in\fq$. Then $f(x)$ is irreducible over $\Fq$ if and only 
if ${\mbox{\Tr}}_{q|2}(\frac{b}{a^2})=1$. 
\end{lemma}

The following lemma is known as Capelli's lemma and can be found in ~\cite{Cohen} too.
\begin{lemma}\label{cohen}
Let $f(x)$ be a degree $n$ irreducible polynomial over $\mathbb{F}_q$, and
let $g(x)\in \mathbb{F}_q[x]$. Then $p(x)=f(g(x))$ is irreducible over 
$\mathbb{F}_q$ if and only if for some root $\alpha$ of $f(x)$ in 
$\mathbb{F}_{q^n}$, $g(x)-\alpha$ is an irreducible polynomial over 
$\mathbb{F}_{q^n}$.
\end{lemma}

Now suppose that both $f(x)=x^2+ax+b$ and $ff(x)$ are irreducible over $\Fq$. We 
show that $fff(x)$ cannot be an irreducible polynomial over $\Fq$. 
Suppose for the contrary that $fff(x)$ is irreducible over $\Fq$. Using 
Lemma~\ref{cohen}, $fff(x)$ is irreducible over $\Fq$if and only if 
$h(x)=f(x)-\alpha$ is irreducible over $\F_{q^4}$ for some root $\alpha$ 
of $ff(x)$ in $\F_{q^4}$. Applying Lemma~\ref{irre-trace}, $h(x)$ is irreducible
over $\F_{q^4}$ if and only if $\Tr_{q^4|2}(\frac{b-\alpha}{a^2})=1$. 
Using the properties of the trace map in one hand  
\begin{equation}\label{Eq:trace}
\Tr_{q^4|2}\left(\frac{b-\alpha}{a^2}\right)=
\Tr_{q^4|2}\left(\frac{b}{a^2}\right)-\Tr_{q^4|2}\left(\frac{\alpha}{a^2}\right)
\end{equation}
and on the other hand
$$
\Tr_{q^4|2}\left(\frac{b}{a^2}\right)=4\Tr_{q|2}\left(\frac{b}{a^2}\right)=0.
$$
Thus from above and Equation~\ref{Eq:trace}, we conclude that 
$\Tr_{q^4|2}\left(\frac{\alpha}{a^2}\right)=1$. But
\begin{eqnarray*}
\Tr_{q^4|2}\left(\frac{\alpha}{a^2}\right)&=
\Tr_{q|2}\left(\Tr_{q^4|q}\left(\frac{\alpha}{a^2}\right)\right)\\
&=\Tr_{q|2}\left(\frac{\Tr_{q^4|q}\left(\alpha\right)}{a^2}\right).
\end{eqnarray*}  
Now suppose that $ff(x)=x^4+c_3x^3+c_2x^2+c_1x+c_0$. Then it is easy to see that $c_3=0$.
Finally, since $\alpha$ is a root of $ff(x)$ in $\F_{q^4}$, we deduce that
 $\Tr_{q^4|q}\left(\alpha\right)=c_3=0$, and hence 
$\Tr_{q^4|2}\left(\frac{\alpha}{a^2}\right)=0$ which is a contradiction.

\section{The case of function fields of characteristic two}
Let $\F_2(t)$ be the rational function field in $t$ over $\F_2$, where $t$ is transcendental 
over $\F_2$, and let $\F_2(t)[x]$ be the polynomial ring over $\F_2(t)$ in the variable $x$.

One expects that a theorem similar to Theorem~\ref{main} to hold for quadratic polynomials
in $\F_2(t)[x]$ but this is not true as the following example shows.

Suppose $f(x)=x^2+t\in\F_2(t)[x]$. Then it is easy to see that
$$
f^{(n)}(x)=x^{2^n}+t^{2^{n-1}}+t^{2^{n-2}}+\cdots+t^2+t.
$$
Now from Eisenstein's criterion for function fields~\cite[Proposition~III.1.14]{Sti} it is easy 
see that for every $n\ge 1$, $f^{(n)}(x)$ is irreducible  over $\F_2(t)$ and hence $f(x)$ is a 
stable polynomial in $\F_2(t)[x]$. 

\section{Acknowledgements}

The author would like to thank Igor Shparlinski for helpful comments on an earlier version of
this note.

\end{document}